\documentclass[10pt]{article}
\usepackage{amsmath,amsthm,amscd,amssymb,mathrsfs,setspace}
\usepackage{latexsym,epsf,epsfig}
\usepackage{epstopdf}
\usepackage{hyperref}
\usepackage[hmargin=3.5cm,vmargin=3.5cm]{geometry}
\usepackage{color}
\newcommand{\ds}{\displaystyle}

\newcommand{\reals}{\mathbb{R}}
\newcommand{\realstwo}{\mathbb{R}^2}
\newcommand{\realsthree}{\mathbb{R}^3}
\newcommand{\xb}{{\bf{x}}}

\newcommand{\pd}{{\partial}}

\newcommand{\Dn}{\partial_{\nu}}

\newcommand{\R}{\mathbb{R}}

\newcommand{\Om}{{\Omega}}

\newcommand{\la}{{\lambda}}

\newcommand{\cH}{\mathcal{H}}

\newcommand{\di}{{\rm div\, }}
\newcommand{\g}{{\nabla}}

\theoremstyle{plain}
\newtheorem{theorem}{Theorem}[section]
\newtheorem{condition}{Condition}

\newtheorem{corollary}[theorem]{Corollary}

\theoremstyle{remark}
\newtheorem{remark}{Remark}[section]
\numberwithin{equation}{section}
\numberwithin{theorem}{section}
\numberwithin{remark}{section}

\title{Mathematical Aeroelasticity: A Survey}
\date{\today}

 \author{\normalsize \begin{tabular}[t]{c@{\extracolsep{.6em}}c@{\extracolsep{.6em}}c@{\extracolsep{.6em}}c}
         Igor Chueshov  & Earl H. Dowell  & Irena Lasiecka & Justin T. Webster  \\
\it        Kharkov National Univ. & \it Duke University & \it Univ. of Memphis    &\it  College of Charleston \\
\it        Kharkov, Ukraine & \it Durham, NC & \it Memphis, TN &\it Charleston, SC\\
 & & \it  IBS, Polish Academy of Sciences  &  \\
\it       chueshov@karazin.ua & dowell@duke.edu&  \it lasiecka@memphis.edu&  \it websterj@cofc.edu
\end{tabular}}

\begin{document}
\maketitle

\begin{abstract} {\noindent
 A variety of models describing the interaction between flows and oscillating structures are discussed. 
 The main aim is to analyze conditions under which structural instability ({\em flutter}) induced by a  fluid flow can be suppressed or eliminated. 
  The analysis provided  focuses on  effects  brought about by:  (i) different plate and fluid boundary conditions, (ii) various regimes for flow velocities: subsonic, transonic, or supersonic, (iii) different modeling of the structure which may or may not account for in-plane accelerations (full von Karman system), (iv) viscous effects, (v) an assortment of models related to piston-theoretic model reductions, and (iv) considerations of axial flows (in contrast to so called normal flows). \\
  The discussion below is based on conclusions reached via  a combination of rigorous PDE analysis, numerical computations, and experimental trials.  
}
\end{abstract}

\section{Introduction}

The main goal of this paper is to provide
a  mathematical survey of various PDE models for fluttering plates \cite{bal0,bolotin,dowellnon,dowell1,gibbs,tang} and associated results, directed at an applied math and engineering readership. The focus here is on {\em modeling concerns}, and we provide a presentation and exposition of mathematical results---when available---and their relationship with known numerical and experimental results concerning  flutter. We do not focus on proofs here, but try to give descriptions of results and a general intuition about why they hold. In addition we mention a handful of mathematical models which are of recent interest in aeroelasticity {\em and also} represent fertile ground from the point of view of mathematical analysis.

 The {\em flutter pheonmenon} is of great interest across many fields. Flutter is a sustained, systemic instability which occurs as a feedback coupling between a thin structure and an inviscid fluid when the natural modes of the structure ``couple" with the fluid's dynamic loading. 
When a structure is immersed in a fluid flow certain flow velocities may bring about a bifurcation in the dynamics of the coupled {\em flow-plate} system; at this point stable dynamics may become oscillatory (limit cycle oscillations, or LCO) or even chaotic. 
A static bifurcation may also occur, known as {\em divergence} or buckling. 

The above phenomena can occur in a multitude of applications: buildings and bridges in strong winds, panel and flap structures on vehicles, and in the human respiratory system (snoring and sleep apnea \cite{sleepapnea1}). Recently, flutter resulting from {\em axial flow} (which can be achieved for {\em low flow velocities}) has been studied from the point of view of energy harvesting \cite{energyharvesting}. Flutter considerations are paramount in the supersonic and transonic regime, with the renewed interest in supersonic flight.
From a design point of view, flutter cannot be overlooked, owing to its potentially disastrous effects on the structure due to sustained fatigue or large amplitude response. 

The field of {\em aeroelasticity} is concerned with (i) producing models which describe the flutter phenomenon, (ii) gaining insight into the mechanisms of flow-structure coupling, (iii) predicting the behavior of a flow-structure system based on its configuration, and (iv) determining appropriate control mechanisms and their effect in the prevention or suppression of instability in the flow-structure system. 
Here we consider flow-plate dynamics corresponding to both  {\it subsonic} and {\em supersonic} flow regimes and a wide array of structural boundary conditions and flow-plate coupling conditions. 

Flow-structure models have attracted considerable attention in the past mathematical literature, see, e.g., \cite{bal4,bal3,b-c, LBC96,b-c-1,chuey,springer,jadea12,dowellnon,webster}
and the references therein.
However, the majority of the work (predominantly in the engineering literature) that
has been done on flow-structure interactions has been devoted to numerical and experimental studies, see,
for instance,
\cite{bal0,bolotin,dowellnon,dowell1,HP02}, and also the survey~\cite{Li03}
and the literature cited there.  Many mathematical studies have been based on linear, two dimensional, plate models  with specific geometries, where  the primary goal was to determine the {\em flutter point} (i.e., the flow speed at which flutter occurs)
\cite{bal0,bolotin,dowell1,HP02,Li03}. See also \cite{bal4,bal-shubov,bal-shubov1,shubov1} for the recent 
studies  of linear models with a one dimensional flag-type structure (beams). This line of work has focused primarily on spectral properties of the system, with particular emphasis on identifying aeroelastic {\em eigenmodes} corresponding to the associated Possio integral equation (addressed classically by \cite{tricomi}). We emphasize that these investigations have been {\em linear}, as their primary goal is to predict the flutter phenomenon and isolate  aeroelastic modes. Given the difficulty of modeling coupled PDEs at an interface \cite{cbms,redbook}, theoretical results have been sparse.
 Additionally, flutter is an inherently nonlinear phenomenon; although the flutter point (the flow velocity for which the transition to periodic or chaotic behavior occurs) can be ascertained within the realm of linear theory, predicting the magnitude of the instability requires a nonlinear model of the structure (and potentially for the flow as well) \cite{dowellnon, B, C, C1}. 

The results presented herein demonstrate that flutter models can be studied from an infinite dimensional point of view, and moreover that meaningful statements can be made about the physical mechanisms in flow-structure interactions {\em strictly from the PDE model}.  
The challenges in the analysis involve (i) the mismatch of regularity between two types of dynamics: the flow and the structure  which are coupled in a hybrid way,  
 (ii) the physically required presence of unbounded or ill-defined terms (boundary traces in the coupling conditions), and (iii) intrinsically non-dissipative  generators
of dynamics, even in the linear case. (The latter are associated with potentially chaotic behavior.) 
One of the intriguing aspects of flow-structure dynamics is that the type of instability, whether static (divergence) or dynamic (LCO), depends both on the plate boundary conditions and on the free-stream (or unperturbed) flow velocity. For example, one observes that if a (two-dimensional) panel is simply supported or clamped on the leading {\em and} trailing edge it undergoes divergence in subsonic flow but flutters in supersonic flow; conversely, a cantilevered panel clamped at one end and free along the others flutters in subsonic flow and may undergo divergence in supersonic flows \cite{K3}.

 This paper primarily addresses the interactive dynamics between a nonlinear plate and a surrounding potential\footnote{For a brief discussion of non-potential and viscous flows we refer to
 Sections \ref{full-vK} and \ref{visc-flow}.}
  flow \cite{bolotin,dowellnon}.
The description of  physical phenomena such as flutter and divergence  will translate into mathematical questions  related to 
 existence of  nonlinear semigroups representing  a  given dynamical system,  asymptotic stability of trajectories, and convergence  to equilibria or to compact attracting sets. 
   Interestingly enough, different model configurations lead to an array of  diverse mathematical issues that involve not only classical  PDEs, but  subtle questions in 
  non-smooth elliptic theory, harmonic analysis, and singular operator theory. For more details concerning  the mathematical theory developed for the flutter models discussed below, see \cite{CDLW1}.

\medskip\par\noindent
{\bf  Notation:}
 For the remainder of the text denote $\xb$ for $(x,y,z) \in \realsthree_+$ or $(x,y) \in \Omega \subset \partial\realsthree_+$, as dictated by context. Norms $||\cdot||$ are taken to be $L_2(D)$ for the domain $D$. 
 The symbol $\nu$ will be used to denote the unit normal vector to a given domain, again, dictated by context. Inner products in $L_2(\realsthree_+)$ are written $(\cdot,\cdot)$, while inner products in $L_2(\R^2\equiv\pd\R^3_+)$ are written $\langle\cdot,\cdot\rangle$. Also, $ H^s(D)$ will denote the Sobolev space of order $s$, defined on a domain $D$, and $H^s_0(D)$ denotes the closure of $C_0^{\infty}(D)$ in the $H^s(D)$ norm denoted by $\|\cdot\|_{H^s(D)}$ or $\|\cdot\|_{s,D}$.  
 We make use of the standard notation for the boundary trace of functions, e.g., for $\phi \in H^1(\realsthree_+)$, $ \gamma (\phi) \equiv tr[\phi]=\phi \big|_{z=0}$ is the trace of $\phi$ on the plane $\{\xb:z=0\}$.

\section{PDE Description of A Fluid-Plate Model}

 The principal model under consideration involves the interaction of a plate with a flow above it. To describe the behavior of the gas (inviscid fluid), the theory of  potential flows \cite{bolotin,dowellnon,dowell1} is utilized. The dynamics of the plate are governed by plate equations with the  von Karman (vK) nonlinearity \cite{springer,ciarlet, lagnese}.
  This model is appropriate for plate dynamics with ``large" displacements (characteristic displacements on the order of a few plate thicknesses),  and therefore applicable for the flexible structures of interest \cite{bolotin,ciarlet,dowellnon,dowell1,lagnese}. There will be configurations in the sequel where vK theory is not valid (e.g., in a flag-like, axial flow configuration), but this will be indicated. 
 Including nonlinearity in the structure of the model is critical, not only for the sake of model accuracy, but also because nonlinear effects play a major mathematical role in controlling ground states; most of  the mathematical results reported in this paper are no longer valid for linearized structures.

The environment considered is $\realsthree_+ = \{(x,y,z) \in\realsthree\, :\, z\ge 0\}$. The thin plate has negligible thickness in the $z$-direction. The unperturbed flow field moves in the $x$-direction at the fixed velocity $U$. The physical constants have been normalized (corresponding to the linear potential flow) such that $U=1$ corresponds to Mach 1 (to simplify the models and to emphasize the mathematical properties of interest).  The plate
 is modeled by  a bounded domain $\Omega \subset \partial \realsthree_+$ with smooth boundary $\partial \Omega = \Gamma$. 
  
 \subsection{Plate} The scalar function $u: \Om\times \R_+ \to \reals$ represents the displacement of the plate in the $z$-direction at the point $(x,y)$ at the moment $t$. The plate is taken with general boundary conditions $BC(u)$, which will be specified later. The dynamics are thus given by:
\begin{equation}\label{plate*}\begin{cases}
       (1-\alpha\Delta)u_{tt}+\Delta_{x,y} ^2 u+D_0(u_t)+f(u)=p(\xb,t) ~\text{in}~\Omega\\
       BC(u) ~ \text{on} ~\partial \Omega,
  \end{cases} \end{equation}
with appropriate initial data $u_0=u(0),~u_1=u_t(0)$. 
The quantity $p(\xb,t)$ corresponds to the aerodynamic pressure of the flow on the top of the plate and, in the standard configuration, is given in terms of the flow: \begin{equation}\label{coupling}\ds p(\xb,t)=p_0(\xb)+(\partial_t +U\partial_x)tr[\phi]\big|_{\Omega},~~\xb\in \Omega.\end{equation}   The quantity $p_0$ represents static pressure applied on the surface of the plate. The function $D_0$ abstractly represents structural damping or a feedback controller. A parameter
$\alpha \ge 0$ represents rotational inertia in the filaments of the plate, and is proportional to the square of the thickness of the plate \cite{lagnese}. The primary case of interest, $\alpha =0$, is said to be {\em non-rotational}; the mathematical effects of $\alpha$ are discussed below.
The vK nonlinearity is given by:  \newline $$\ds f(u)=-[u, v(u)+F_0],
$$ where $F_0$ is a given {\em in-plane load} (of sufficient regularity). 
      The vK bracket  above corresponds to
$$
\ds
[u,w] = u_{xx}w_{yy}+
u_{yy}w_{xx} -
2 u_{xy}w_{xy},
$$ 
and
the Airy stress function $v(u_1,u_2) $ solves the  elliptic
problem
\begin{equation*}
\Delta^2 v(u_1,u_2)+[u_1,u_2] =0 ~~{\rm in}~~  \Omega,\quad \Dn v(u_1,u_2) = v(u_1,u_2) =0 ~~{\rm on}~~  \Gamma,
\end{equation*}
 (the notation $v(u)=v(u,u)$ is employed).
 \par 
 The following types of boundary conditions for the displacement $u$ are typical for flow-structure interaction models.
 \begin{itemize}
\item[(C)] {\it  Clamped boundary conditions} (corresponding to a panel element) take the form 
\begin{equation*} u = \Dn u = 0 ~\text{ on } \Gamma. \end{equation*} 
In this case, the plate is  considered to be embedded  in or affixed to  a large rigid body.
\item[(FC)]  Let the boundary be partitioned in two pieces: $\Gamma_0$ and $\Gamma_1$. {\it  Free-clamped boundary conditions} (possibly corresponding to a flap, flag, or cantilevered airfoil) are given by
 \begin{equation*}
\begin{cases} u=\partial_{\nu}=0~~\mathrm{on}~~ \Gamma_0\\ \Delta u + (1-\mu) B_1 u = D_1(\Dn u_t)~~\mathrm{on}~~ \Gamma_1,
 \\
\partial_{\nu} \Delta u + (1-\mu) B_2 u - \mu_1 u = D_2(u,u_t)  + \alpha \Dn u_{tt} ~~\mathrm{on}~~ \Gamma_1, \end{cases} \end{equation*}
 where the boundary operators $B_1$ and $B_2$
are given by \cite{lagnese}:
\begin{equation}
\begin{array}{c}
B_1u = 2 \nu_1\nu_2 u_{xy} - \nu_1^2 u_{yy} - \nu_2^2 u_{xx}=-\partial_{\tau}^2u-\nabla \cdot \nu(\xb)\Dn u\;, \\
\\ 
B_2u = \partial_{\tau} \left[ \left( \nu_1^2 -\nu_2^2\right) u_{xy} + \nu_1
\nu_2 \left( u_{yy} - u_{xx}\right)\right]\,=\partial_{\tau}\partial_{\nu}\partial_{\tau}u;%
\end{array}
\end{equation}
  The parameter $\mu_1$ is nonnegative and $0<\mu<1$ is the Poisson modulus. The operators $D_i$ for $i=1,2$ are possible energy damping/dissipation functions. 
\end{itemize}
 \subsection{Flow}
 The scalar function $\phi: \realsthree_+\times \reals \to \reals$ is the linear {\em flow potential} and satisfies:
\begin{equation} \label{flow}\begin{cases}
     (\partial_t+U\partial_x)^2\phi=\Delta_{x,y,z} \phi ~\text{in} ~\realsthree_+\\
    FC(\phi),
        \end{cases} \end{equation}
       with appropriate initial data $\phi_0$ and $\phi_1$\footnote{Note that, as mentioned above, $U$ is normalized by the speed of sound and thus is what engineers refer to as the {\em Mach number.}}. The term $FC(\phi)$ represents the {\em flow boundary conditions} or {\em interface conditions}. We will consider two primary flow conditions:
        \begin{itemize}
        \item[(NC)] The standard flow boundary conditions---the {\em full Neumann condition}---are henceforth denoted (NC). This is typically utilized when a majority of the plate boundary is clamped or hinged. It takes the form:
        \begin{equation}\label{NC}
        \partial_z\phi \big|_{z=0} = (\partial_t+U\partial_x)u_{\text{ext}},~~\text{on }~ \mathbb{R}^2=\partial \mathbb{R}^3_+,
        \end{equation}
         where the subscript `ext' indicates the extension by zero for functions defined on $\Omega$ to all of $\mathbb R^2$. 
        \item[(KJC)] The {\em Kutta-Joukowsky} flow condition, henceforth denoted (KJC), is a mixed (Zaremba-type \cite{savare}) condition which will be discussed  below.  We introduce the {\em acceleration potential} $\psi = (\partial_t+U\partial_x)\phi$ for the flow potential $\phi$. Then the  {\it Kutta-Joukowsky} flow condition on $\realstwo$ is
        \begin{equation}\label{KJC}
        \partial_z \phi = (\partial_t+U\partial_x)u ~ \text{ on } ~\Omega;~~ \psi = 0  \text{ on } ~ \realstwo \setminus \Omega.
        \end{equation}
        \end{itemize}
   In both cases above the interface coupling on the surface of the plate occurs in a Neumann type boundary condition known as the {\em downwash} of the flow. 
   The choice of flow conditions is itself dependent upon the plate boundary conditions, and this is determined by application. 
 \section{Panel Flutter  Model}
 Our primary interest here is the following PDE system:
\begin{equation}\label{flowplate}\begin{cases}
(1-\alpha\Delta)u_{tt}+\Delta^2u+ku_t +f(u)= p_0+\big(\partial_t+U\partial_x\big)tr[\phi] & \text { in } \Omega\times (0,T),\\
u(0)=u_0;~~u_t(0)=u_1,\\
u=\Dn u = 0 & \text{ on } \partial\Omega\times (0,T),\\
(\partial_t+U\partial_x)^2\phi=\Delta \phi & \text { in } \realsthree_+ \times (0,T),\\
\phi(0)=\phi_0;~~\phi_t(0)=\phi_1 & \text { in } \realsthree_+\\
\Dn \phi = -\big[(\partial_t+U\partial_x)u (\xb)\big]_{\text{ext}}& \text{ on } \realstwo_{\{(x,y)\}} \times (0,T).
\end{cases}
\end{equation}
We focus on this panel flutter model as a starting point for the analysis of coupled flow-plate systems. After providing results for this model (a majority of extant mathematical results hold for this model), we describe extensions and recent work to accommodate other configurations in the sequel. 
  \subsection{Functional Setup and Energies}
 \noindent Finite energy constraints (as dictated by the physics of the model) manifest themselves in the natural  topological requirements on the solutions $\phi$ and $u$ to \eqref{flowplate}. Letting $\cH$ correspond to the state space for plate displacements $u(t)$ (depending on boundary conditions), and considering the space $L^{\alpha}_2$ with  $||\cdot||_{L_2^{\alpha}(\Omega)}^2 = \alpha ||\nabla \cdot||^2+||\cdot||^2$ for plate velocities and thus $$L_2^{\alpha} (\Omega)  \equiv 
 H^1(\Omega),~\text{for}~\alpha>0,~ \text{and} ~ L_2(\Omega) ,~\text{for}~\alpha =0.$$ Solutions we consider here will have the properties: $u \in C(0,T; \cH)\cap C^1(0,T;L_2^{\alpha}(\Omega)); ~ \phi \in C(0,T;H^1(\realsthree_+))\cap C^1(0,T;L^2(\realsthree_+)).$   To set up the model in a dynamical systems framework, the principal state space is taken to be $$Y = Y_{fl}\times Y_{pl} \equiv \big(W_1(\realsthree_+) \times L_2(\realsthree_+)\big)\times\big( \cH \times L^{\alpha}_2(\Omega)\big)$$
 where $W_1(\R^3_+)$ will denote the homogeneous Sobolev space of order $1$. We will also consider a stronger space: $$Y^s \equiv H^1(\realsthree_+)\times L_2(\realsthree_+) \times \mathcal H \times L_2^{\alpha}(\Omega).$$
 
  \begin{remark}[Parameters]The central parameters in the PDE analysis of the general models above are the rotational inertia parameter $\alpha \ge 0$ and the flow velocity $U$. First, $\alpha$ is critical in that the model of primary physical interest in the theory of large deflections takes $\alpha=0$, yet the essential nonlinearity {\em and} interactive flow terms are {\em not of a compact nature} in this case. Indeed, the presence of rotational terms provides an additional regularizing effect 
 on the transverse velocity of the plate, $u_t$ \cite{b-c,b-c-1,springer};  this, in turn, leads to several desirable mathematical properties such as compactness of trajectories, gain of derivatives,  etc.
   It is also clear that, with respect to $U$, when $U>1$, there is a loss of spatial ellipticity of the principal part of the flow operator $\partial_t^2-\Delta+U^2\partial_x^2$, and near $U=1$ there is full $x$ degeneracy in the model. \end{remark}
   
\begin{remark} For the remainder of the text the case of clamped plate boundary conditions (C) with standard Neumann flow conditions (NC) is referred to as the {\bf standard panel} configuration.
\end{remark}
The subsonic case, $U \in [0,1)$, is initially considered. The following ``energies" are observed with the standard {\em velocity multipliers} (we suppress the dependence on $t$ in the expressions below):
\begin{align}
E_{pl}(u) = & ~\dfrac{1}{2}[||u_t||_{L_2^{\alpha}(\Omega)}^2+ ||\Delta u ||_{0,\Omega}^2+\frac{1}{2}||\Delta v(u)||_{0,\Omega}^2]-\langle F_0,[u,u]\rangle_{\Omega}+\langle p,u\rangle_{\Omega} \\
E_{fl}(\phi)=  &~\dfrac{1}{2}[||\phi_t||^2_{0,\realsthree_+}+||\nabla \phi||_{0,\realsthree_+}^2-U^2||\phi_x||^2_{0,\realsthree_+}]\\
E_{int}(u,\phi) =&  ~U\langle\phi,u_x\rangle_{\Omega};\hskip1cm\mathcal E = E_{pl}+E_{fl}+E_{int}
\end{align}
where $E_{int}$ represents the (not necessarily positive) {\em interactive} ``energy".

 In the discussion below, we will encounter strong (classical), generalized (mild), and weak (variational) solutions.
In obtaining existence and uniqueness of solutions, semigroup theory is utilized; this necessitates the use of \textit{generalized} solutions. These are strong limits of strong solutions and  satisfy an integral formulation of (\ref{flowplate}) (they are called \textit{mild} by some authors). We note that generalized solutions are also weak solutions and satisfy a variational formulation of the flow-plate system (see, e.g., \cite[Section 6.5.5]{springer} and \cite{webster}). For the exact definitions of 
 strong and generalized  solutions
we refer  to \cite{springer,jadea12,supersonic,webster}. 

\subsection{Well-Posedness of Nonlinear Panel Model}\label{wellp}
 
\subsubsection{Overview of Previous Results}
Well-posedness results in the past literature deal mainly with the dynamics possessing some regularizing effects. This has been accomplished by either accounting for non-negligible rotational inertia effects \cite{b-c,b-c-1,springer} and strong damping of the form $-\alpha \Delta u_t$, or by incorporating helpful thermal effects  into the structural model \cite{ryz,ryz2}.  In the cases listed above, the natural structure of the dynamics dictate that the plate velocity has the property $u_t \in H^1(\Omega)$,
 which provides  the needed regularity 
for the applicability of many standard tools in nonlinear analysis. One is still faced with the low regularity of boundary traces,
due to the failure of the Lopatinski conditions \cite{sakamoto}.

The first contribution to the mathematical analysis of the problem  is \cite{LBC96,b-c-1} (see also \cite[Section 6.6]{springer}), where the case  $\alpha > 0$ (rotational) is fully treated.  The method employed in \cite{LBC96,b-c-1,springer}
relies on sharp microlocal estimates for the solution to the wave equation driven by $H^{1/2}(\Omega) $ Neumann  boundary data given by~ $u_t + U u_x$. This gives~ $\phi_t|_{\Omega} \in
L_2(0,T; H^{-1/2}(\Omega))$  \cite{miyatake}.  Along with an explicit solver for the three dimensional flow equation (the Kirchhoff formula), and a Galerkin approximation for the structure, one may construct  a fixed point for the appropriate
 solution map. This method  works well for all values of $U\ne1$. 

  Related ideas were used  more recently  in the thermoelastic flow-plate interactions in \cite{ryz,ryz2}; in this case the dynamics enjoy
$H^1(\Omega)$ regularity of the velocity $u_t$ (independent of $\alpha$) due to the analytic regularity  induced by thermoelasticity \cite{redbook}. 
When $\alpha =0$, and we take no additional smoothing in the model, the  estimates corresponding to the above approaches become singular, destroying the  applicability of the  previous techniques.  The  main  mathematical difficulty of this problem is  the presence of  the   boundary  terms:  $ (\phi_t +
 U \phi_x)|_{\Omega} $  ~acting as the aerodynamic pressure on the plate. When $U =0$, the corresponding  boundary terms exhibit monotone behavior with respect to   the natural energy  inner product (see \cite[Section 6.2]{springer} and \cite{cbms}), which is topologically equivalent  to  the topology of the state space $Y$.   The latter  enables  the use of monotone operator theory \cite{springer} (and references therein).  However, when
  $U > 0$ this is no longer true, and other techniques must be employed to obtain solutions \cite{jadea12,webster}.  

In \cite{bal0,bal4,shubov3,shubov1} (and references therein)  a 
linear {\em airfoil}  immersed in a subsonic flow is considered. The wing is taken to have 
a high aspect ratio thereby allowing for the suppression of the span variable, 
and reducing the analysis to individual chords normal to the span. By reducing the problem to a one dimensional analysis, many technical hangups are avoided, and Fourier-Laplace analysis is  effective. The analysis of  two dimensional structures requires very different tools. 

\subsubsection{Existence and Uniqueness Results}
We begin with  an overview of the well-posedness results for the the standard panel model which simultaneously covers both subsonic and supersonic cases (and thus is a weaker result than what appears in the sequel).

	\begin{theorem}\label{th:supersonic}
Consider  problem in (\ref{flowplate}) with $U\ne1$, $p \in L_2(\Omega), F_0 \in H^{3 +\delta}(\Omega)$.  Let $ T > 0 $ and
\begin{equation}\label{space-Y}
( \phi_0, \phi_1; u_0, u_1 ) \in  
Y = Y_{fl}\times Y_{pl} \equiv \big(W_1(\R^3_+) \times L_2(\realsthree_+)\big)\times\big( H^2_0(\Omega) \times L^{\alpha}_2(\Omega)\big)
\end{equation}
 where 
 $L_2^{\alpha} (\Omega) =  
 H^1_0(\Omega)$ if $\alpha>0$ and  $L_2^{\alpha} (\Omega)=L_2(\Omega)$ when $\alpha =0$.
Then there exists unique generalized solution
\begin{equation}\label{phi-u0reg}
(\phi (t), \phi_t(t); u(t), u_t(t)) \in C([0, T ],  Y).
\end{equation}
This solution is also weak and generates
  a  nonlinear  continuous  semigroup
$S_t : Y \rightarrow Y$.
\end{theorem}
The proof of the above theorem in the general (supersonic) case makes use of a change of state variable: the dynamic variable $\psi = \phi_t+U\phi_x$ is considered as the second flow state. In this case the adapted energies are:
\begin{align}
\widehat E_{fl}= ~ \dfrac{1}{2}[||\nabla \phi||_{0,\realsthree_+}^2]+||\psi||^2_{0,\realsthree_+},~~
\widehat E_{int} = & ~0,~~
\widehat {\mathcal E} =  ~E_{pl}+\widehat E_{fl}+\widehat E_{int}.
\end{align}
Generalized solutions satisfy the identity is given by:
\begin{equation}\label{enr}\ds
\widehat{\mathcal E}(t) + \int_s^t\int_{\Omega} ku_t^2 + \textcolor{red}{U}\int_0^t\langle u_x,\psi\rangle_{\Omega} d\tau = \widehat {\mathcal E}(s).\end{equation}
\begin{remark}
Note that when $U > 1$ the equation (\ref{flowplate}) displays degeneracy  of ellipticity  of the stationary problem. This results in the loss of dissipativity in the energy relation  (\ref{enr}). 
 In addition, the boundary term $\psi|_{\Omega} $ is not well defined on the energy space. 
 Handling of these  obstacles requires the development of an appropriate ``hidden" trace regularity theory for the aeroelastic potential $\psi$,
 which plays pivotal role in the proof of Theorem \ref{th:supersonic}.\end{remark}
	
	In the subsonic case,  instead, there is no spatial ``degeneracy" in the flow equation, leading to a stronger well-posedness result:
\begin{theorem}\label{th:subsonic} Suppose $U<1$, $p \in L_2(\Omega)$ and $F_0 \in H^{3+\delta}(\Omega)$.
 Assume 
that initial data sytisfy \eqref{space-Y}. Then there exists unique generalized solution to \eqref{flowplate}
\begin{equation}
(\phi (t), \phi_t(t); u(t), u_t(t)) \in C([0, T ],  Y).
\end{equation}
This solution is also weak and generates
  a  nonlinear  continuous  semigroup
$S_t : Y \rightarrow Y$.
 Every generalized solution is also weak.
Moreover these solutions to \eqref{flowplate} generates
a  nonlinear continuous semigroup $(S_t,Y)$. 
Any generalized (and hence weak) solution to \eqref{flowplate},  satisfies the global-in-time bound \begin{equation}\label{energybound}
\sup_{t \ge 0} \left\{((1-\alpha \Delta)u_t, u_t)_{0,\Omega}^2+\|\Delta u\|_{\Omega}^2+\|\phi_t\|_{\realsthree_+}^2+\|\nabla \phi\|_{\realsthree_+}^2 \right\}  < + \infty.\end{equation}
and thus $(S_t,Y)$ is uniformly bounded in time, i.e., there exists a constant C such that $$ ||S_t (y_0)||_Y \leq C (||y_0||_Y), ~~t>0.$$
The dynamics satisfy the energy identity 
\begin{equation}\label{ener-rel-subs}
{\mathcal{E}}(t)+  \int_s^t \int_{\Omega} k u_t^2={\mathcal{E}}(s),~~ t\ge s,
\end{equation}
\end{theorem}

\begin{remark} This theorem remains true if  the von Karman nonlinearity is replaced by any nonlinearity $f: H^2(\Omega) \cap H_0^1(\Omega) \to L_2(\Omega)$ which is (a) locally Lipschitz, and (b) the nonlinear trajectories possess a uniform bound in time which typically holds when nonlinearity satisfies some coercivity properties, see \cite{ch-l}.
\end{remark}

The proof of Theorem \ref{th:subsonic} given in 
\cite{springer} and
\cite{webster} relies on three  principal ingredients: (i) renormalization of the state space which yields shifted dissipativity for the dynamics operator (which is nondissipative in the standard norm on the state space); (ii) the sharp regularity of Airy's stress function, which converts a supercritical nonlinearity into a critical one  \cite{springer} (and references therein): 
\begin{equation}\label{airy}
||v(u)||_{W^{2,\infty}(\Omega) } \leq C ||u||^2_{H^2(\Omega)};
\end{equation}
and (iii) control of low frequencies for the system by the nonlinear source, represented by the inequality
\begin{equation}
||u||^2_{L_2(\Omega)} \leq \epsilon [ ||u||^2_{H^2(\Omega)^2} + ||\Delta v(u)||^2_{L_2(\Omega) } ] + C_{\epsilon} .
\end{equation}
An alternative proof based on a {\em viscosity method} can be found in \cite{jadea12}. 
\medskip\par

  In comparing the results obtained for supersonic and subsonic cases,  there are two major differences at the qualitative level, in addition to an obvious fact that $E_{fl}$ is no longer positive when $U > 1$:
{\it First}, the regularity of strong solutions obtained in the subsonic case \cite{jadea12,webster} coincides with regularity expected for classical solutions. In the supersonic case, there is a  loss of differentiability in the flow  in the tangential  $x$ direction, which then propagates to the  loss of differentiability in the structural variable $u$. As a consequence, strong solutions do not exhibit sufficient regularity in order to perform the  needed calculations.
  To cope with the problem a special regularization procedure was introduced in \cite{supersonic},
  where  strong solutions are approximated by sufficiently regular functions, though not solutions to the given PDE.
  The limit passage allows one to obtain the needed estimates valid for the original solutions \cite{supersonic}.
{\em Secondly}, in the subsonic case, one shows that the solutions are {\it bounded} in time,
  see \cite[Proposition 6.5.7]{springer} and also \cite{jadea12,webster}.
  This property  cannot be shown in the supersonic analysis.  The leak of the energy in the energy relation can not be compensated by the nonlinear terms. 

\subsection{Long-time Behavior and Attracting Sets}
 It is known \cite{dowell1} (and references therein) that, although the flow coupling may introduce potential instability\footnote{Calculations performed at the physical level demonstrate that the effect of the flow can be destabilizing in certain regimes---see \eqref{negdamp} below.}, it can also help bound the plate dynamics and may assist in the dissipation of plate energy associated to high frequencies \cite{C,C1}. The work in \cite{delay} is an attempt to address these observations rigorously, via the  result in Theorem \ref{rewrite} below, which reduces the full flow-plate dynamics to a plate equation with a delay term. This result is implemented and the problem is cast within the realm of PDE with delay (and hence our state space for the plate contains a delay component which encapsulates the flow contribution to the dynamics). 
One would like to show that such dynamics, as described above, are stable in the sense that trajectories converge asymptotically to a ``nice" set. Unlike parabolic dynamics, there is no a priori reason to expect that hyperbolic-like dynamics can be asymptotically reduced to a truly finite dimensional dynamics. By showing that the PDE dynamics converges to a finite dimensional (compact) \textit{attractor} it effectively allows the asymptotic-in-time reduction of the analysis of the infinite dimensional, unstable model  to a finite dimensional set upon which classical control theory can be applied. We note that owing to the non-dissipative terms present in flutter models, global attractors associated to plate's flutter dynamics (as well as the set of equilibria points for the full flow-plate system) can have non-trivial structure. See \cite{dowellA} and references cited therein, as well as \cite{dowellB,dowellC,dowellD} for the engineering discussion of these properties.
\begin{remark}
Based on numerical results, it is well known that for a panel whose width is much larger than its chord, that four to six structural modes are sufficient to give very good accuracy to predict flutter and limit cycle oscillations. Indeed, for transonic and subsonic flows even one or two modes may be sufficient. As the width of the panel is decreased more and more modes are required, however. Also if the plate is under tension, then the number of modes needed increases, and in the singular limit of infinite tension the number of modes required also goes to infinity \cite{B,dowellA}.\end{remark}

The main  point  of the treatment in \cite{delay} is to demonstrate the existence and finite-dimensiona\-lity of a global attractor for the reduced plate dynamics (see Theorem \ref{rewrite}) {\em in the absence of any imposed damping mechanism} and in the absence of any regularizing mechanism generated by the dynamics itself.
Imposing no mechanical damping, taking $\alpha =0$, and assuming the flow data are compactly supported we obtain the primary theorem in \cite{delay}:
\begin{theorem}\label{th:main2}
Suppose $0\le U \ne 1$, $\alpha =0$, $F_0 \in H^{3+\delta}(\Omega)$ and $p_0 \in L_2(\Omega)$.
 Then there exists a compact set $\mathscr{U} \subset Y_{pl}$ of finite fractal dimension such that $$\displaystyle \lim_{t\to\infty} d_{Y_{pl}} \big( (u(t),u_t(t)),\mathscr U\big)
 \displaystyle=\lim_{t \to \infty}\inf_{(\nu_0,\nu_1) \in \mathscr U} \big( ||u(t)-\nu_0||_2^2+||u_t(t)-\nu_1||^2\big)=0$$
for any weak solution $(u,u_t;\phi,\phi_t)$ to \eqref{flowplate} with
initial data
$$
(u_0, u_1;\phi_0,\phi_1) \in Y,~~Y\equiv H_0^2(\Omega)\times L_2(\Omega)\times W_1(\realsthree_+)\times L_2(\realsthree_+),
$$
which are
localized  in $\R_+^3$ (i.e., $\phi_0(\xb)=\phi_1(\xb)=0$ for $|\xb|>R$ for some $R>0$). Additionally, there is the extra regularity $\mathscr{U} \subset \big(H^4(\Omega)\cap H_0^2(\Omega)\big) \times H^2(\Omega)$.
\end{theorem}

A key to obtaining attracting sets is the representation of the flow on the structure via a delay potential (see Section 3.3 in \cite{springer}).  Reducing this full flow-plate problem to a delayed von Karman plate is the primary motivation for our main result and permits a starting-point for long-time behavior analysis of the flow-plate system, which is considerably more difficult otherwise. The exact statement of this key reduction is given in the following assertion:

\begin{theorem}\label{rewrite}
Let the hypotheses of Theorem~\ref{th:main2} be in force, and $(u_0,u_1;\phi_0,\phi_1) \in \cH \times L_2(\Omega) \times H^1(\realsthree_+) \times L_2(\realsthree_+)$. Assume that there exists an $R$ such that $\phi_0(\xb) = \phi_1(\xb)=0$ for $|\xb|>R$.  Then the there exists a time $t^{\#}(R,U,\Omega) > 0$ such that for all $t>t^{\#}$ the weak solution $u(t)$ to (\ref{flowplate}) also satisfies the following equation:
\begin{equation}\label{reducedplate}
u_{tt}+\Delta^2u-[u,v(u)+F_0]=p_0-(\partial_t+U\partial_x)u-q^u(t)
\end{equation}
with
\begin{equation}\label{potential}
q^u(t)=\dfrac{1}{2\pi}\int_0^{t^*}ds\int_0^{2\pi}d\theta [M^2_{\theta} u_{\text{ext}}](x-(U+\sin \theta)s,y-s\cos \theta, t-s).
\end{equation}
Here, $M_{\theta} = \sin\theta\partial_x+\cos \theta \partial_y$ and \begin{equation}\label{delay} t^*=\inf \{ t~:~\xb(U,\theta, s) \notin \Omega \text{ for all } \xb \in \Omega, ~\theta \in [0,2\pi], \text{ and } s>t\}
\end{equation} with $\xb(U,\theta,s) = (x-(U+\sin \theta)s,y-s\cos\theta) \subset \realstwo$.
\end{theorem}

Thus, after some time, the behavior of the flow can be captured by the aerodynamical pressure term $p(t)$ in the form of a delayed forcing.  Theorem~\ref{rewrite} allows us to suppress the dependence of the problem on the flow variable $\phi$. The flow state variables $(\phi,\phi_t)$  manifest themselves in our rewritten system via the delayed  character of the problem; they appear  in the initial data for the delayed  component of the plate, namely $u^t\big|_{(-t^*,0)}$. 
Here we emphasize that the structure of aerodynamical pressure posited in the hypotheses leads to the velocity term  $-u_t$.
We may utilize this as natural damping appearing in the structure of the reduced flow pressure. 
\par
 The reduced model displays the following features: (i) it does not have gradient structure (due to dispersive and delay terms), (ii) the delay term appears at the critical level of regularity. However,  despite of the lack of gradient structure, compensated compactness methods allow to ``harvest" some compactness properties from the reduction, so the ultimate dynamics does admit global attracting set. 

The proof of Theorem \ref{th:main2} requires modern tools and new long-time behavior technologies applied within this delay framework. Specifically, the approach mentioned above (and utilized in \cite{springer,ryz,ryz2}) {\em does not apply} in this case.  A relatively new technique \cite{ch-l,springer} allows one to address the asymptotic compactness property for the associated dynamical system without making reference to any gradient structure of the dynamics (not available in this model, owing to the dispersive flow term). In addition, extra regularity and finite dimensionality of the attractor is demonstrated via a quasi-stability approach  \cite{springer}
(see also \cite{ch-dqsds} for recent  developments of quasi-stability techniques).
The criticality of the nonlinearity and the lack of gradient structure prevents one from using a powerful technique of {\it backward smoothness of trajectories}  \cite{springer}, where smoothness is obtained near the equilibria points and propagated forward.  Without a gradient structure for this model, the associated attractor may have complicated structure (not being strictly characterized by the equilibria points).

In the presence of additional  damping imposed on the structure, it is reasonable to expect that the entire evolution (both plate and flow) is strongly stable in the sense of strong convergence to  equilibrium states. Such results are shown for the model  above (in the same references) only for the case of {\em subsonic flows}, $U<1$. The key to obtaining such results is a {\em viable energy relation} and a priori bounds on solutions which yield {\em finiteness of the dissipation integral}.
 When $U>1$, recent calculations reported in \cite{jfs} indicate that such results are false, even for low supersonic speeds.
The following result provides a stabilization of the full subsonic flow-plate dynamics (as opposed to the result above we obtain for all flow velocities, but only the structural dynamics). 

For this we recall
  the global bound in \eqref{energybound}
which along with the energy identity in \eqref{ener-rel-subs} 
 allows us to obtain the following: \begin{corollary}\label{dissint}
Let $ k> 0 $, $0\le U<1$  and $\alpha=0$. Then the dissipation integral is finite. Namely, for a generalized solution to \eqref{flowplate} we have $$  \int_0^{\infty} \|u_t(t)\|_{0,\Omega}^2 dt \leq  
K_{u,k} < \infty.$$
\end{corollary}
\noindent From this finiteness, we can show the following result \cite{conequil1} for {\em smooth initial data}. Let $\mathcal N$ denote the set of stationary solutions to \eqref{flowplate} (for the 
 existence and properties see \cite{springer}).
\begin{theorem}\label{regresult}
Let $0\le U<1$ and $\alpha=0$.  Assume $p_0 \in L_2(\Omega)$ and $F_0 \in H^{4}(\Omega)$. Assume that $\mathcal N$ is {\em discrete}\footnote{This is generically true; see \cite{springer,conequil1} for a more detailed discussion.} Then for all $ k>0$, any  solution $(u(t),u_t(t);\phi(t),\phi_t(t))$ to the flow-plate system \eqref{flowplate} with 
initial data
$$
(u_0,u_1;\phi_0,\phi_1) \in (H_0^2\cap H^4)(\Omega) \times H_0^2(\Omega)\times H^2(\mathbb R^3_+) \times H^1(\mathbb R^3_+).
$$
that are
 spatially localized in the flow component (i.e., there exists a $\rho_0>0$ so that for $|\xb|\ge \rho_0$ we have $\phi_0(\xb) = \phi_1(\xb)=0$) has the property that 
 \begin{align*}\lim_{t \to \infty}\left\{\|u(t)-\hat u\|^2_{H^2(\Omega)}+\|u_t(t)\|^2_{L_2(\Omega)}+\|\phi(t)-\hat \phi\|_{H^1( K_{\rho} )}^2+\|\phi_t(t)\|^2_{L_2( K_{\rho} )} \right\}=0
 \end{align*} 
 for  some $(\hat u,\hat \phi) \in \mathcal N$, and 
 for any   $\rho>0$, where 
 $K_\rho=\{ \xb\in\R_+^3\, : \|\xb\|\le \rho\}.$
\end{theorem}
The above result remains true for finite energy initial data (only in $Y$) with convergence in a {\em weak sense}. 
Additionally, for the system with rotational inertia present $\alpha>0$ (and corresponding strong damping) \cite{chuey,springer} {\em or} with thermal effects present \cite{ryz,ryz2} the analogous {\em strong} stabilization result holds  for finite energy initial data 
$$
(u_0,u_1;\phi_0,\phi_1) \in Y=H_0^2(\Omega) \times L^{\alpha}_2(\Omega)\times W_1(\mathbb R^3_+) \times L_2(\mathbb R^3_+).
$$

The fact that $ k>0 $ is needed in order to obtain convergence to equilibria is corroborated by the counterexample  \cite{igor}, 
 which shows that, with $ k =0$, periodic solutions may remain in the limiting dynamics. In fact, taking $ k >0 $ allows one to prove that the entire flow-structure system is a gradient-type system. More specifically, the simplified example in \cite{igor} considers the  linear plate model without any damping $k =0$, coupled to the flow via matching velocities.  (Thus, in this example, $U$ is taken to be zero.) The boundary conditions imposed on then plate are simply supported (hinged). In this scenario,  periodic solutions may persist, and in \cite{igor} they are explicitly constructed for a specific tubular domain. Whether the same result  holds for clamped plate boundary data---is an open question. However, this example indicates the necessity of introducing mechanical damping in the plate model {\em if one expects a strong convergence to equilibria} of the full flow-structure system.

 In view of the above  result---Theorem \ref{regresult}---with smooth initial data, we see that {\em any damping} imposed on the structure seems to eliminate the flutter in the subsonic case. (This is consistent with experiment for clamped plates, as discussed above.) From a physical point of view, the aerodynamic damping for a plate in subsonic flow is much lower than in supersonic flow, so the plate tends to oscillate in a neutrally stable state in the absence of aerodynamic damping or structural damping \cite{dowellA}.

We also note that
when the vK nonlinearity is replaced by Berger's nonlinearity, and the damping coefficient $ k > 0$ is taken sufficiently large, it has been shown \cite{conequil2} that 
all {\em finite energy} trajectories converge strongly to equilibria, and hence the flutter is eliminated. Though we have Theorem \ref{regresult} for {\em smooth initial data}, the analogous result for the vK nonlinearity   is  still unknown, unless an additional static damping controller is used \cite{conequil1}. Indeed, in \cite{conequil1} we consider finite energy data and the vK or Berger nonlinearity; when large static {\em and} large viscous damping is active in the model, we can show strong convergence of the entire flow-plate system to equilibria.

\section{Piston Theory and the Hypersonic Limit}\label{piston}
A {\em piston-theoretic} approximation/simplification of the dynamics arises in considering the surface of the plate being acted upon by ``piston". The piston-theoretic term is a nonlinear function of the {\em downwash} term (in the case of interest here is the Neumann condition for the flow on appearing on the RHS of the plate equation). This term is ``expanded" via a series, and the traditional piston theory (taken to be valid for $U>\sqrt{2}$) considers only the linear term in the expansion. In some early piston-theoretic literature \cite{BA62}, polynomial terms in ~$-[u_t+Uu_x]$~ are retained up to the third order.
Thus, standard, linear piston theory (or law of plane sections \cite{bolotin,dowellnon,oldpiston}) replaces the {\em acceleration potential of the flow},~ $\psi= tr[\phi_t+U\phi_x]$ with ~$-[u_t+Uu_x]$ in the RHS of the reduced plate equation \eqref{reducedplate}.

The following overview is given in \cite{venedeev2}: ``Piston theory initially was obtained as an asymptotic expansion of the exact expression [for the flow contribution] as $U \to \infty$. In later studies it was shown that piston theory is valid starting from $U \approx 1.7$...Coupled mode (coalescence) [classical] flutter is observed for values of $U$ beyond this critical value, whereas single mode flutter can occur in the range $U \in (1,1.7)$."

Formally, we could  arrive at the aforementioned, standard (linear) piston theory model by
 utilizing the reduction result in Theorem \ref{rewrite}; a bound  exists \cite[p.334]{springer}: \begin{equation}
||q^u(t)||_{0,\Omega}^2 \le \dfrac{C}{U} \int _{t-\frac{c}{U}}^t||\Delta u(\tau)||_{\Omega}^2 d\tau
\end{equation}
for $U$ sufficiently greater than 1. The constants $c,C>0$ depend only on the diameter of $\Omega$. This indicates that the delay term, $q^u$ decays (in some sense) as $U$ increases. Hence, 
it is reasonable to guess that $q^u$ can be neglected in the case of large speeds $U$.\footnote{See also \eqref{negdamp} for another (physically motivated) expression for 
the ``piston'' term.} 
So we arrive at the following model
\begin{equation}\label{plate-stand}\begin{cases}
u_{tt}+\Delta^2u+f(u) = p^*(\xb,t) ~~ \text { in } ~\Omega\times (0,T), \\
u=\Dn u =0  ~~\text{ on } ~ \partial\Omega\times (0,T),
\end{cases}
\end{equation}
where ~~$\ds p^*(\xb,t)=p_0-[u_t+Uu_x]$~~ is the piston-theoretic dynamic pressure.

  This model was intensively studied in 
the literature for different types of boundary conditions and frictional damping forces, see \cite{springer} and also \cite{ch-l} for related second order abstract models. The typical result is the following assertion:
\begin{theorem}
Under the conditions imposed on the model
the equations \eqref{plate-stand} generates a dynamical system in the state space $H_0^2(\Omega)\times L_2(\Omega)$, possessing a compact global attractor of finite fractal dimension which is also more regular than the state space, in the sense of a gain of two derivatives in each component.	
\end{theorem}
Existence of compact and finite dimensional attractor for the structure yields an important reduction of an infinite dimensional, hyperbolic-like PDE system to a finite dimensional dynamics (ODE) system in the asymptotic limit. Such a result  can  not  to be (a priori) predicted by a specific or numerical experiment. In addition, the rigorous mathematical validation of this fact is far from obvious, since the reduced plate  system {\em is not a gradient dynamical systems}, and there is no natural compactness built in into the dynamics. In fact, one of the key roles is played by the nonlinearity, whose presence is critical for   proving an existence of the absorbing set. 

 \begin{remark} Free plate boundary conditions have been recently studied in \cite{bociu} in the context of piston theory with boundary (and/or geometrically constrained) dissipation.\end{remark}

A key avenue of research for supersonic flow-plate models in the standard panel configuration \cite{supersonic} (see Section \ref{wellp}) is an explicit comparison between solutions to the piston-theoretic plate model and solutions to the clamped flow-plate interaction, as $U$ becomes sufficiently large. 
More specifically: one should show how that the associated solutions to \eqref{reducedplate} and \eqref{plate-stand} discussed above coincide for (i) $t \to \infty$ with $U$ fixed and sufficiently large, or (ii) as $U \to \infty$ (the so called hypersonic limit) on some {\em arbitrary time interval $[T_1,T_2]$} (where $T_2$ is perhaps $\infty$).\footnote{
The reference \cite{hypersonic} contains a result for the behavior of solutions to the piston-theoretic plate 
(with RHS given by $p^*$ above) as $U\to \infty$. The result, however, is only valid for arbitrarily small time intervals, and hence does not provide information about the behavior of solutions for arbitrary $t$.}

\begin{remark}\label{nonpist} We note that, as mentioned above, some classical piston-theoretic analyses \cite{BA62} retain polynomial terms in ~$-[u_t+Uu_x]$~ up to the third order. Adding third powers of ~$-[u_t + U u_x]$~  will  lead  to a model with nonlinear, monotone (cubic-type) damping in $u_t$
(see also some discussion in \cite[p.161]{ch-l}). Such damping has already been considered in wave and plate models \cite{springer}, and is known to induce stronger stability properties in mitigating effects of nonconservative terms in the equation.  It would be interesting to study this situation in the context of the piston model under consideration.  \end{remark}

 In addition to the classical, piston-theoretic model discussed above (in \eqref{plate-stand} and Remark \ref{nonpist}), there are other {\em piston theories}.
For instance, if we a priori  assume that we deal with ``low'' frequencies regimes for the plate, then we can use the following expression for the aerodynamic pressure (see, e.g., \cite{B,dowellA} and the references therein):
\begin{equation}\label{negdamp}
p^{**}(\xb,t) =p_0- \dfrac{U}{\sqrt{U^2-1}}\left(\left[\dfrac{U^2-2}{U^2-1}\right]u_t+Uu_x\right).
\end{equation}
in the supersonic case ($U>1$).
We note that this term has a fundamentally different structure, and its contribution to the dynamics is markedly different than the RHS of \eqref{plate-stand}.
It is commonly recognized that the damping due to the aerodynamic flow is usually substantially larger in {\em magnitude} than the damping due to the plate structure. However, as we can see,  the aerodynamic damping can be negative as well as positive.
 For instance, in the case $1<U < \sqrt{2}$
we obtain negative damping \cite{B,dowellA}. This  has a definite
impact on the dynamics, producing instabilities in the system. We also note that for very large $U$  ($U \rightarrow \infty$) the model in \eqref{negdamp} coincides with the standard  one  described by (\ref{plate-stand}).

In a recent paper \cite{newpiston}, the piston-theoretic approach is re-examined utilizing an additional term in the asymptotic expansion in the inverse Mach number. The result is an aeroelastic pressure on the surface of the plate of the form 
\begin{equation}\label{pistoon}
p^{***}(\xb,t)=p_0(\xb)-(\partial_t+U\partial_x)u(t)+\frac{1}{U}\int_0^x(\partial_{tt}+U\partial_{tx})u(\xi,t) d\xi
\end{equation}

The reference \cite{newpiston} discusses this model in the context of an extensible beam and two dimensional flow, but indicates that the model may be generalized. Additionally, the derivation shown in \cite{newpiston} is valid for $U$ values only slightly larger than $U=1$. This may provide a key for studying low supersonic flow-plate interactions. From a mathematical point of view, one may also notice that the formula in (\ref{pistoon}) indicates a cancellation between the stabilizing part of the flow, represented by $u_t$, and the additional integral term.
\begin{remark} The force $p^*$ in equation \eqref{plate-stand} is the result of expanding the full linear potential flow theory in powers of frequency for fixed $U$ (i.e., Mach number), while equation (\ref{pistoon}) is the result of expanding the full potential flow theory in in powers of $1/U$ (or inverse Mach number). These contribute very different effects to the dynamics, unless higher order terms are retained in frequency and $1/U$.\end{remark}

\section{Recent Configuration of Interest}\label{recent}

 A configuration which supports multiple models of recent interest is now described. The principal component of this configuration is the existence of large portion of the plate boundary $\Gamma_2\equiv \partial \Omega_2$ which is free. 
 \begin{remark} As a general remark, we note that changes in the configuration (e.g., plate boundary conditions) can have immense impacts on the overall character of dynamics. This mathematical observation is confirmed experimentally \cite{dowellE} as well. \end{remark}
 Taking a free-clamped plate boundary (as given by (FC)) allows one to consider (in generality) the situation of a {\em cantilevered wing}, as well as a {\em mostly clamped panel}, or a flap/flag in so called {\em axial flow}. When one considers a free plate coupled to fluid flow (unlike the clamped case) the key modeling issues correspond to (i) the validity of the structural nonlinearity \cite{D} and (ii) the aerodynamic theory near the free plate boundary (which may exhibit in-plane displacement).

\begin{center}
 \includegraphics[scale=.22]{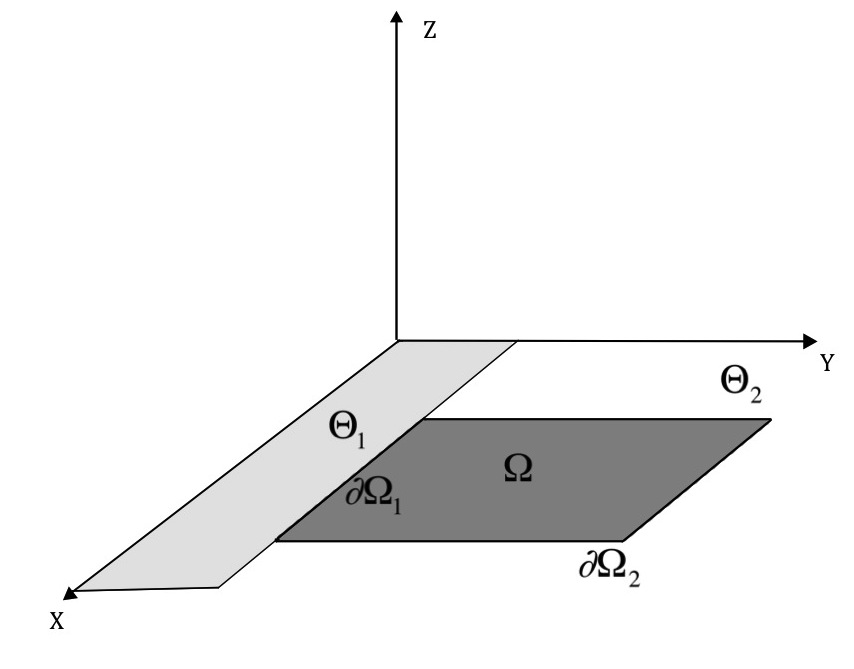}
\end{center}

In this configuration, a natural {\em flow} interface boundary condition arises and is called {\em Kutta-Joukowsky flow condition} (KJC), as described above and in \cite{bal0,K1,dowellA,K2,K3}.  Simply put, the (KJC) corresponds to a zero pressure jump off the wing and at the trailing edge \cite{bal4,dowell1}. The (KJC) has been implemented in numerical aeroelasticity as mechanism for  ``...removal of a velocity singularity at some distinguished point on a body in unsteady flow" \cite{K1}. 
From an engineering standpoint, the (KJC) is required to provide a unique solution for the potential flow model for a lifting surface, and gives results in correspondence with experiment, i.e. the pressure difference across the trailing edge is zero and the pressure difference at the leading edge is a near maximum. Studies of viscous flow models in the limit of very high Reynolds numbers lend support to the (KJC).

\begin{remark}\label{axial} The configuration above arises in the study of airfoils. In this case, we refer to {\em normal flow} (along the $x$-axis). Another  related configuration referred to as {\em axial flow}, which takes the flow occurring along $y$ axis. Physically, the orientation of the flow can have a dramatic effect on the occurrence and magnitude of the oscillations associated with the flow-structure coupling. This manifests itself, specifically, in the choice of nonlinearity modeling the plate (or possibly beam) equation. \end{remark}

\subsection{Kutta-Joukowsky Flow Condition and Free Plate Boundary Conditions}\label{KJ}
\noindent The Kutta-Joukowsky conditions (KJC) described above \eqref{KJC} are {\em dynamic and mixed} in nature.  
These mixed type flow boundary conditions are taken to be accurate for plates in the clamped-free configuration. 
The  configuration below represents an attempt to model oscillations of a plate which is {\em mostly free}. The dynamic nature of the flow conditions corresponds to the fact that the interaction of the plate and flow is no longer static at the free edge.
In this case we take the free-clamped plate boundary conditions, and the mixed flow boundary conditions:
\begin{equation}\label{physical}\begin{cases}
u=\Dn u = 0,  &\text{ on }  \pd\Om_1\times (0,T) \\
{B}_1 u =0,   {B}_2 u =0,~ &\text{on}~ \pd\Om_2\times (0,T)\\
\Dn \phi = - (\partial_t+U\partial_x)u, ~& \text{on}~\Omega\times (0,T)\\
\Dn \phi= 0, & ~\text{on}~\Theta_1\times(0,T)\\
\psi=\phi_t+U\phi_x=0, &~ \text{on}~\Theta_2\times(0,T)\end{cases}
\end{equation}
where $\pd\Om_i,~ i =1,2$ are complementary parts of the boundary $\partial \Omega $ , and ${B}_1$, ${B}_2$ represent moments and shear forces, given by \cite{lagnese} (and earlier here) and $\Theta_i$ extend in the natural way into 
the remainder of the $x-y$ plane)

The implementation of free-clamped plate boundary conditions is extremely important in the 
modeling of airfoils. However, treating coupled fluid-plate problems that involve a free plate is technically challenging. This is due to the  loss of sufficient  regularity of the boundary data imposed for the flow  (the failure of the Lopatinski condition). Clamped boundary conditions assumed on the boundary of the plate allow for smooth extensions to $\R^2$ of the (NC) conditions satisfied by the flow. In the absence of these, one needs to approximate the original dynamics in order to construct sufficiently smooth functions amenable to PDE calculations.  Preliminary calculations indicate that (as the physics dictate) free plate boundary conditions are in fact {\em more compatible} with the (KJC) flow conditions. 

\begin{remark}
With regard to (i) above, we note that when we are in the case of {\em normal flow} (as opposed to {\em axial flow} (as in Remark \ref{axial}) the scalar von Karman model is largely still viable \cite{dowell1}.
Additionally, any configuration where the free portion of the (FC) plate boundary conditions is ``small" with respect to the clamped portion will satisfy the hypotheses for the theory of large deflections, and hence, von Karman theory is applicable. 
\end{remark}

 \subsection{An Exploratory Case: (KJC) in Subsonic Panel Configuration}
 Arguably, the  Kutta-Joukowsky  boundary conditions for the flow (\ref{KJC})  are the {\em most important} when modeling an airfoil  immersed in a flow of gas \cite{bal,bal0}. 
Not surprisingly, these boundary conditions are also  the most challenging  from mathematical   stand point. This was recently brought to the fore in an extensive treatise \cite{bal0}. 
Various aspects of the problem---in both subsonic and supersonic regime---have been  discussed  in \cite{bal0} in the context of  mostly one dimensional structures. 

The aim of this section is to address the mathematical problem posed by (KJC), by putting them within the framework of modern harmonic analysis. We consider the ``toy" problem of (KJC) conditions coupled with clamped plate boundary conditions below. Indeed, the recently studies  \cite{fereisel,KJ} show how the flow condition (KJC) interacts with the clamped plate in subsonic flows {\em in order to develop a suitable abstract theory} for this particular  flow condition. 
  The resulting papers \cite{fereisel,KJ}  give well-posedness of this fluid-structure interaction configuration.
   Though the analysis is subsonic for (KJC), utilizing the flow energy from the supersonic panel (as above in Theorem \ref{th:supersonic}) is effective in the abstract setup of this problem.  
In fact,  even in the subsonic case, the analysis of semigroup generation proceeds  through the technicalities developed  earlier for supersonic case \cite{supersonic}. The key distinction from the analysis of the clamped  flow-plate interaction (owing to the dynamic nature of the boundary conditions) is that a {\em dynamic flow-to-Neumann map} must be utilized to incorporate the boundary conditions into the abstract setup. The regularity properties of this map are critical in admitting techniques from abstract boundary control \cite{redbook}, and are determined from the Zaremba elliptic problem \cite{savare}.
  The necessary trace regularity hinges upon the invertibility of an operator which is analogous (in two dimensions) to the finite Hilbert transform \cite{bal0,bal4,tricomi}.   And this is a critical additional element of the challenging  harmonic analysis brought about  by the (KJC). 
When the problem is reduced in dimensionality to a beam structural model, this property can be demonstrated and our analysis  has parallels  with that in \cite{bal0,bal4} (and older references which are featured therein). Specifically, one must invert the finite Hilbert transform in $L_p(\Omega)$ for $p \ne 2$; in higher dimensions, this brings about nontrivial (open) problems in harmonic analysis and the theory of singular integrals.

 From the mathematical point of view
the difficulty lies, again,  at the level of the  linear theory. In order to deal with   the effects of the {\it unbounded}
traces $tr[\psi] $ in the energy relation microlocal calculus is necessary. This
has been successfully accomplished in \cite{supersonic} where clamped boundary conditions in the supersonic case
were considered. However, in the case of (KJC)  there  is an additional difficulty that involves
``invertibility" of  {\it finite} Hilbert (resp. Riesz) transforms.  This latter property is known to fail within $L_2$ framework, thus it is necessary to build the  $L_p$ theory, $p \ne 2 $.  This was for the first time observed in \cite{bal4}
and successfully resolved in the  one dimensional case. However, any progress to higher dimensions depends
on the validity of the corresponding harmonic analysis result developed for finite Riesz transforms.

  The full system with  (\ref{KJC}) can be  written in terms of aeroelastic potential variable as:
\begin{equation}\label{flowplate2}\begin{cases}
u_{tt}+\Delta^2u+f(u)= p_0+tr[\psi] & \text { in } \Om\times (0,T),\\
u(0)=u_0;~~u_t(0)=u_1,\\
u=\Dn u=0 & \text{ on } \pd\Om\times (0,T),\\
(\partial_t+U\partial_x)\phi=\psi & \text { in } \realsthree_+ \times (0,T),\\
(\partial_t+U\partial_x)\psi=\Delta \phi & \text { in } \realsthree_+ \times (0,T),\\
\Dn \phi = -(\partial_t+U\partial_x)u & \text{ on } \Omega   \times (0,T),\\
(\partial_t+U\partial_x)  \phi = 0 & \text{ on } \realstwo \backslash \Omega \times (0,T),\end{cases}
\end{equation}
with the given initial data $(u_0,u_1,\phi_0, \psi_0 ) $. 
 We make use of the flow acceleration multiplier $(\partial_t+U\partial_x)\phi \equiv \psi$ for the {\em subsonic flow}, taken with (KJC). 
Thus for the flow dynamics, instead of $(\phi;\phi_t)$ we again have
the state variables  $(\phi;\psi)$ (as in the supersonic panel case)   which then leads to a {\it  non dissipative}  energy balance.

 Our result is formulated under the following regularity condition (to be discussed later).
 \begin{condition}[{\bf Flow Trace Regularity}]\label{le:FTR0}
  Assume that  $\phi(\xb,t)$  satisfies \eqref{flow} and the Kutta-Joukowsky condition (KJC),
   and $\partial_ttr[\phi],~ \partial_xtr[\phi]  \in L_2(0,T;H^{-1/2-\epsilon}(\realstwo))~~~~\forall\, T>0.
$
In addition,  $\psi = \phi_t + U \phi_x $ satisfies the estimate
\begin{equation}\label{trace-reg-est-M0}
\int_0^T\|tr[\psi](t)\|^2_{H^{-1/2 -\epsilon} (\realstwo)}dt\le C_T\left(
E_{fl}(0)+
 \int_0^T\| \Dn \phi(t) \|_{\Omega}^2dt\right)
\end{equation}
\end{condition}
We note that the regularity required by Condition \ref{le:FTR0} is exactly the one that will make energy relation meaningful. 
To wit: $u_x(t) \in H^1_0(\Omega)  $ and $ tr[\psi] \in L_2(0,T; H^{- \theta } ), ~\theta > 1/2 $ defines the correct  duality pairing in the tangential direction on $\mathbb R^2$. 
It is also at this point where we use the fact that $u$ satisfies  $ u=0  $ on $\partial \Omega $ . Thus simply supported and clamped boundary conditions imposed on the structure  fully cooperate with this regularity.  
The principal result of \cite{KJ} (see also \cite{fereisel}) reads as follows: 
\begin{theorem}\label{T1}
With reference to the model (\ref{flowplate2}), with  $0 \leq U < 1$:
  Assuming the trace regularity Condition \ref{le:FTR0} holds for the aeroelastic potential $\psi = (\phi_t+U\phi_x)$,
  there exists a unique finite  energy solution which exists on any $[0,T]$.  
  This is to say, for any $T > 0 $, 
  $(\phi, \phi_t, u, u_t ) \in C(0, T; Y )$  for all initial data  $(\phi_0, \phi_1, u_0, u_1 ) \in  Y $. 
  This solution depends continuously on the initial data.  
\end{theorem}
The proof of Theorem \ref{T1} follows the technology developed for the supersonic case in \cite{supersonic} and is given in \cite{KJ}. 
In view of Theorem \ref{T1} the final conclusion on generation of nonlinear semigroup is pending upon verification of the flow regularity Condition. 
While the complete solution  to this question is still unavailable, and pending further progress in the theory of finite Riesz transforms, we can provide positive answer in the case when $\Omega$ is one dimensional. 
In fact, this positive assertion also follows  a-posteriori from direct  calculations given in \cite{bal}, and based on the analysis of finite Hilbert transforms. 
The arguments given \cite{KJ} are independent and more general, with potential adaptability   to multidimensional cases. 
The corresponding result is formulated below. 
\begin{theorem}\label{1dtracereg}
Assume $\Omega=(-1,1)$ (suppressing the span variable $y$). Then the flow trace regularity in Condition \ref{le:FTR0} holds for $\phi$.
In this case  the semiflow defined by (\ref{flowplate2}) taken with (KJC) generates a continuous semigroup.
\end{theorem}

\begin{remark}
As discussed above, the generation of semigroups for an arbitrary three dimensional flow is subjected to the validity of the trace Condition \ref{le:FTR0}. While it is believed that this property should be generically true, at the present stage this appears to be an open question in the analysis of singular integrals and depends critically on the geometry of $\Omega$ in two dimensions. 
\end{remark}

\begin{remark}\label{Bal1}
We note that  the  regularity of the aeroelastic potential required by  Condition \ref{le:FTR0} in one dimension, follows from the analysis in \cite{bal0,bal4}, where the author proves that 
aeroelastic potential $ \psi \in L_2(0, T , L_q(\Omega) ) $ for $ q < 4/3 $ . Since  for $p > 4 $ there exists $\epsilon > 0 $ 
such that $$H^{1/2 + \epsilon} (\Omega) 
\subset L_p(\Omega) , ~~p > 4 , ~~\text{dim}~~\Omega \leq 2,$$  and one then obtains  that $L_q(\Omega) \subset H^{-1/2 - \epsilon}(\Omega) $ with $q < 4/3 $. \end{remark}
\begin{remark}
The loss of $1/2$ derivative in the characteristic region was already observed and used in the analysis of regularity of 
the aeroelastic potential for the Neumann problem with supersonic velocities $U$ \cite{supersonic}. 
However, in the case of (KJC)  there is an additional loss, due to the necessity of inverting finite Hilbert transform which forces to work with $L_p$ theory for $ p < 2 $. This is due to the fact that  {\it finite} Hilbert transform  is invertible on $L_p$, $p < 2$, rather than $ p =2$. 
\end{remark}

\subsection{Full von Karman Model}\label{full-vK}
\noindent In the considerations above (as well as most in mathematical aeroelasticity) either linear plate theory or scalar von Karman equations have been implemented to model the large deflections of the plate in the flow-plate interaction. One possible mechanism for greatly improving the generality in modeling is to utilize the so called {\em full von Karman equations} which account for in-plane accelerations as well as transverse accelerations of the plate. In the case of the full von Karman equations, the in-plane displacement acts a system of elasticity, and the structure of the nonlinearity, though coupled, behaves in a {\em local fashion} unlike the scalar version of von Karman. However, considering this system greatly complicates certain aspects of PDE analysis. 

\begin{align}\label{fullkarman}
(1-\alpha\Delta)u_{tt}+\Delta^2u-&\text{div}\big(\mathcal C[\varepsilon(\mathbf w)+f(\nabla  u)]\nabla  u \big) =p _1(\xb,t)\\
\mathbf{w}_{tt}-&\text{div} \big(\mathcal C[\varepsilon(\mathbf w)+f(\nabla u)]  \big) = P_2(\xb,t)\\
BC(u,\mathbf{w})& ~~\text{ on } ~ \Gamma;~~~f(s)=~\frac{1}{2}s \otimes s,
\end{align}
where $u$ is the transversal displacement of the plate and $w=(w^1;w^2)$ is the tangential (in-plane) displacements.
The term $p_1(\xb,t)$ again incorporates flow coupling in the transverse sense (the aerodynamic pressure on the surface of the plate), and the vector $P_2$ allows for flow coupling in the in-plane sense as well. The term $\varepsilon$ corresponds to the standard strain tensor given by $\varepsilon(w) = \frac{1}{2}[\nabla w+\nabla w^T]$, and the fourth order symmetric tensor $\mathcal C$ is given by 
$$\mathcal C(\varepsilon) \equiv \dfrac{E}{(1-2\mu)(1+\mu)}[\mu ~\text{trace}( \varepsilon I)+(1-2\mu)\varepsilon],$$ where $\mu$ has the meaning of the Poisson modulus, and $E$ Young's modulus.
\par 
In the case of full vK system we need also to take into account compatibility of tangential (in-plane) movements of the plate with the corresponding dynamics of the gas flow.
This leads to the following model for the viscous flow with states $(u,w,p,v)$ (for more detail, see \cite{fereisel}):
\begin{align}
 &
   p_t+ U p_{x}+\di\, v=0\quad~~ {\rm in}~~ \R^3_+
   \times\R_+,\label{den-eq1U}
   \\[2mm]
 &
   v_t+ U v_{x} -\mu_f\Delta v -(\mu_f+\la_f)\g \di\, v + \g p =0~~ \quad {\rm in}~ \R^3_+
   \times\R_+,\label{flu-eq1U}
\end{align}
where $p$ is the pressure and $v=(v^1;v^2;v^3)$ 
is a small perturbation
of the gas velocity field $(U;0;0)$
Here $\mu_f$ and $\la_f$ are (non-negative) viscosity coefficients (which  vanish
in the case of  invisid fluid).
We need also supply \eqref{den-eq1U} and \eqref{flu-eq1U} with appropriate boundary conditions.  We can choose non-slip boundary conditions, for instance.
They  \textcolor{red}{are given by } 
\begin{equation}\label{non-slip}
(v,\tau)=
\begin{cases}
0 ~&~~ {\rm on}~\R^3_+\setminus \Om, \\
(w_t,\tau)	~& ~~{\rm on}~\Om,
\end{cases}
~~\mbox{and}~~~ v^3=
\begin{cases}
0 ~& ~~{\rm on}~\R^3_+\setminus \Om, \\
u_t+Uu_{x}	~ &~~ {\rm on}~\Om,
\end{cases}
\end{equation}
where $\nu$ is the outer normal, $\tau$ is a tangent direction (see also some discussion 
 in \cite{dowell1} where the explanation concerning the term ~$u_t+Uu_{x}$~
on the boundary is given).
Aerodynamical impact in the full von Karman system \eqref{fullkarman} is described by the forces
\[
p_1(\xb,t)= T^{33}=2\mu_f v^3_{z}+\la_f\, \di v-p
\]
and
\[
P_2(\xb,t)=(T^{13};T^{23})= \mu_f\left(v^1_{z}+v^3_{x}; v^2_{z}+v^3_{y}\right),
\]
where  $T=\{T^{ij}\}_{i,j=1}^3$ is the stress tensor of the fluid.
We need also equip the equations above with 
initial data for all variables $u,w,v,p$.
\par
In the case of 
\textbf{invisid  compressible} fluid ($\mu_f=0$ and  $\la_f=0$) we can 
introduce  (perturbed) velocity potential\footnote{As in \cite{Chu2013-invisid}
we can also deal with equations for velocities 
directly without introducing potential.}
$(v=-\g \phi$) which satisfies the equation  
\[
(\partial_t+U\partial_x)^2\phi=\Delta_{x,y,z} \phi ~~\text{in} ~\realsthree_+\\
\]    
with the boundary conditions
\begin{equation}\label{non-slip-phi}
\frac{\pd\phi}{\pd\tau}=
\begin{cases}
0 ~& {\rm on}~\R^3_+\setminus \Om, \\
-(w_t,\tau)	~& {\rm on}~\Om,
\end{cases}
~~\mbox{and}~~~ \frac{\pd\phi}{\pd\nu}=
\begin{cases}
0 ~& {\rm on}~\R^3_+\setminus \Om, \\
-(u_t+Uu_{x})	~ & {\rm on}~\Om,
\end{cases}
\end{equation}
Moreover, in this case, the pressure has the form 
$p= u_t+Uu_x$ and thus  
\[
p_1(\xb,t)= \big(\partial_t+U\partial_x\big)tr[\phi] 
~~\mbox{and}~~
P_2(\xb,t)= 0.
\]
To the best of our knowledge the fluid-structure  
model in the setting  described above was not 
studied yet at mathematical level. 
We only note that
over the past 25 years there have been a handful of treatments which address the PDE theory of the full von Karman plate system. We do not here provide an extensive overview of well-posedness and stability analyses of the full von Karman equation; the references \cite{koch,fullkarman3,fullkarman2} address: well-posedness of classical solutions, well-posedness of energy type solutions for $\alpha >0$, and stability in the presence of thermal effects or boundary dissipation (respectively). At present, well-posedness of energy-type solutions for $\alpha=0$ is an open problem. Recently, the analyses in \cite{cr-full-karman,berlin11} successfully demonstrated well-posedness and the existence of attractors for certain fluid-structure models involving the full von Karman equations
(in the case $U=0$). 
The  fluid-plate models under consideration in these references make use of the full von Karman plate model coupled with a three dimensional incompressible fluid in a bounded domain (this corresponds the case when equation \eqref{den-eq1U} is neglected and equation \eqref{flu-eq1U} is taken 
with $U=0$ and $\di v=0$). 
An elastic membrane bounds a portion of the fluid-filled cavity. This model is motivated by blood flowing through a large artery, and in the modeling of the liquid sloshing phenomenon \cite{slosh2} (fluid in a flexible tank). It assumes a homogeneous, viscous, incompressible fluid modeled by Stokes flow in \cite{cr-full-karman}. In \cite{cr-full-karman},  well-posedness of energy-type solutions was shown along with the existence of a compact attractor in the presence of frictional damping. The damping is critical (Ball's method \cite{jadea12,supersonic}) for obtaining the attractor. 

\subsection{Viscous flows}\label{visc-flow}
In some situations, a general model 
with a velocity fluid field of the general form (like in \eqref{den-eq1U} and \eqref{flu-eq1U}
with the nonzero viscosity coefficients $\mu_f$ and $\la_f$) may arise. 
This is a case when it necessary to take into account viscosity effects near the oscillating plate, which can be important for transonic flows, and possibly for hypersonic flows

Neglecting in-plane displacements in the \eqref{fullkarman} model of the previous subsection, we arrive at the system:
\begin{equation}\label{flowplate-compr}\begin{cases}
(1-\alpha\Delta)u_{tt}+\Delta^2u+ku_t +f(u)= tr_\Omega \big[2\mu_f v^3_{z}+\la_f\, \di v-p\big]
 & \text { in } \Omega\times (0,T),\\
u(0)=u_0;~~u_t(0)=u_1  & \text { in } \Omega,\\
u=\Dn u = 0 & \text{ on } \partial\Omega\times (0,T),\\
p_t+ U p_{x}+\di\, v=0 & {\rm in}~~ \R^3_+
   \times (0,T),
   \\
   v_t+ U v_{x} -\mu_f\Delta v -(\mu_f+\la_f)\g \di\, v + \g p =0 & {\rm in}~ \R^3_+
   \times (0,T) \\
v(0)=v_0;~~p(0)=p_0 & \text { in } \realsthree_+
\\
v =\left(0;0; \big[(\partial_t+U\partial_x)u (\xb)\big]_{\text{ext}}\right)& \text{ on } \realstwo_{\{(x,y)\}} \times (0,T).
\end{cases}
\end{equation}
The mathematical theory for this model is not yet well-developed.
Some mathematical results on well-posedness and long-time dynamics are available for this case when a fluid fills a bounded or tube type domain and $U=0$. See \cite{Chu2013-comp} and references therein.
For the corresponding model in the incompressible case ($\di v=0$) and with $U=0$
 we refer to  \cite{ChuRyz2011,ChuRyz2012-pois}.

\subsection{Axial Flow}

\noindent In the case of {\em axial flow}, mentioned in Remark \ref{axial}, a beam or plate is clamped on the {\em leading edge} and free elsewhere is described. To provide a clear picture of the dynamics, consider the figure at the beginning of Section \ref{recent} with the over-body flow occurring in the {\em $y$-direction},  as opposed to being in the direction of the plate's chord---$x$-direction.\footnote{We remark that this convention is typically reversed $x\leftrightarrow y$ in the engineering literature.}
This axial configuration, owing to LCO response to low flow velocities \cite{bal0}, is that which has been considered from the point of view of energy harvesting. This would be accomplished by configuring an axially oriented flap on land or air vehicles with a piezo device which could generate current as it flutters \cite{energyharvesting}.
From \cite{tang} the appropriate structural equations of motion here correspond to those of a thin pipe conveying fluid, and many aspects of such dynamics mirror those we are investigating \cite{D,tang}. Expressions for the kinetic and potential energy of the beam and plate are given in \cite{dowellG,dowellF}. Note that structural nonlinearities occur in both the inertia (kinetic energy) and stiffness (potential energy) terms for the structural equations of motion.

The recent theoretical and experimental work of Tang and Dowell \cite{dowellG,dowellF} has been encouraging. In this work a new nonlinear structural model has been used based upon the inextensibility assumption and the comparison between theory and experiment for the LCO response has been much improved over earlier results \cite{tang}. The study of a linear aerodynamic model, combined with the new nonlinear structural model,  is worthy of more rigorous mathematical attention.

\section{The Transonic Regime}

Much of the technical discussion presented above excludes the sonic velocity $ U =1$; in practice, the vicinity of $U\approx1$ is known as the {\em transonic regime}.
Indeed, for $U=1$ the analysis  provided above for the panel configuration breaks down in the essential way \cite{transnon}.
While numerical work   predicts appearance of shock waves \cite{C}, to our best knowledge  no mathematical treatment of this problem is  available at present.

In the literature \cite{dowell1}, it is noted that the flow model discussed above in \eqref{flow} is not accurate in the regime $U \approx 1$\footnote{The strongly transonic region occurs in the range $U \in [.95,1.2]$.}; indeed, the flow equations require additional specificity when the flow velocity nears the transonic regime. We note that the base model in \eqref{flow} becomes interesting, as a {\em degenerate wave equation} appears when $U=1$. This will clearly produce diminished regularity in the $x$ direction. However, the references \cite{transnon,C} suggest that in the transonic regime a fully nonlinear fluid must be considered. From \cite{transnon}: ``...neglecting either the structural or the fluid nonlinearities can lead to completely erroneous stability predictions." The reference \cite{dowellA} indicates that in the case of the standard panel configuration the appropriate flow equation---due to the {\em local Mach number} effect---for transonic dynamics is of the form:
\begin{equation}
\left(\partial_t+\mathbf{1}\partial_x\right)^2\phi=\Delta \phi -\left(\mathbf{1} \right)\phi_{x}\phi_{xx}.
\end{equation}
Here  $U$ is replaced with $\mathbf{1}$ indicating that $U \approx 1 \pm \epsilon$. This nonlinear fluid equation introduces new mathematical challenges which must be addressed. In the case of supersonic or subsonic flows, the definite sign of the spatial flow operator associated to the $x$ variable can be exploited. Additionally, the nature of the dynamics become quasilinear in the flow potential $\phi$. Thus, near the transonic barrier the flow equation becomes $x$-degenerate and quasilinear. 

Experimentally, one notes that near $U=1$ there are many peculiarities, including the possibility of hysteresis \cite{tang} in $U$, and various flow instabilities, including shocks \cite{transnon,C}. Such shocks may actually induce flutter. In fact, recent numerical investigations \cite{jfs,venedeev2} investigate the emergence of {\em single mode flutter} in low supersonic speeds, and the possibility of stabilizing this mode via structural damping. However, the results in \cite{transnon} indicate that shocks emanating from the flow may actually induce more complex coupled mode (coalescence \cite{venedeev2}) flutter. Also in the transonic range, viscous fluid boundary layer effects are apt to be more important \cite{C1,dowellA}.

It is essential that interaction between structural nonlinearities and the aerodynamic nonlinearities be accounted for. In \cite{bendixen}, numerical studies of the interaction between structural and aerodynamic nonlinearities have been investigated. According to this study, so called {\em traveling wave flutter} occurs in the transonic regime. In fact, for a simply supported panel, at $20,000~ft$ of altitude, with a thickness-chord ratio of $0.004$ traveling wave flutter has been recorded  in the form of wave packets; these packets  change  shape and evolve in time during the movement from the leading to the trailing edge. Near the panel there has been evidence of shocks forming, moving in concert with the panel motion, and in some cases {\em inducing} flutter.

From the mathematical point of view, one notes that a key issue is the regularity of flow solutions in the $x$-direction. However, as $\Delta-(\mathbf{1})\partial_x^2$ becomes degenerate, a term providing additional information in the $x$ direction of the flow appears.  Specifically, in energy calculations, the quasilinearity gives rise to a conserved quantity in the flow energy of the form $\ds \Pi(\phi) = \dfrac{-U}{6}||\phi_x||^3_{0,\realsthree_+},$ so 
\begin{equation}\label{transflow}E_{fl}= \frac{1}{2}\left[||\phi_t||^2+||\nabla \phi||^2-(\mathbf{1})^2||\phi_x||^2-\frac{(\mathbf{1})}{3}||\phi_x||^3\right].\end{equation} This quantity, although conserved, has a ``bad" sign, which is a potential indicator for shock phenomena. It seems clear that the analytical techniques will differ for $U = 1-\epsilon$ and $U=1+\epsilon$. Moreover, if well-posedness can be established, owing to physical results, we anticipate that the qualitative properties of these flow (and hence flow-structure) solutions will be quite distinct. 

From the point of view of quasilinear theory, the term $\phi_x \phi_{xx} $ provides dominant control of the operator $D^2_x$, and cannot be neglected in this transonic case. In fact, the corresponding PDE problem can be viewed from the point of view of the Tricomi operator \cite{bal0} changing from elliptic to hyperbolic, depending on the direction of the change in deflection. This is an exciting question to study mathematically; the inroads provided by numerics and experiment here are instrumental guiding forces.

\section{Concluding Remarks}
In conclusion we recall shortly 
 several observations made in \cite{CDLW1}  which, on one hand confirm experimentally and numerically some of the mathematical findings and, on the other hand, raise open questions and indicate  new avenues for mathematical research. For more detailed discussion we refer to our paper \cite{CDLW1}.
\begin{enumerate}
\item{\bf  Stability induced by the flow in reducing dynamics to finite dimensions}:
Rigorous mathematical analysis reveals that the stabilizing effect of the flow reduces the structural dynamics to a {\it  finite dimensional } setting. It is remarkable that such conclusion is obtained directly
from mathematical considerations---as it is not definitively arrived at by either numerics or experiment. 
\item{\bf Non-uniqueness of final nonlinear state}:
  It is possible that there exist several \textit{locally} stable equilibria in the global attractor of the system. This explains why
  the {\em buckled plate} in an aerodynamic flow does not have a final, unique nonlinear state.  
\item
{\bf Surprisingly subtle effects of boundary conditions}: 
As the structure's boundary conditions are changed, so is the dynamic stability/instability of the system. Experimental and numerical and theoretical studies confirm this.
\item
{\bf Limitations of VK theory and new nonlinear plate theory}: 
When the leading edge is clamped and the side edges and trailing edge are free, the vK theory can be no longer accurate.  
However, a novel, improved nonlinear plate theory has been developed and explored computationally and correlated with experiment. This provides a  challenging opportunity for mathematical analysis.
\item
{\bf Aerodynamic theory}:
 It may be helpful to prove, mathematically, some of the results
 related to the justification of the recent developments in piston theory.
    This issue is also important for modeling the pressure due to solar radiation, which has a similar mathematical form to that of piston theory \cite{dowellH,dowellI}. It  has been of recent interest in the context of interplanetary transportation using solar sails. 
\end{enumerate}

\section{Acknowledgment}
The authors would like to dedicate this work to Professor A.V. Balakrishnan, whose pioneering and insightful work on flutter brought together engineers and mathematicians alike. 

E.H. Dowell was partially supported by the National Science Foundation with grant NSF-ECCS-1307778.
I. Lasiecka was partially supported by the National Science Foundation with grant NSF-DMS-0606682 and the United States Air Force Office of Scientific Research with grant AFOSR-FA99550-9-1-0459.
 J.T. Webster was partially supported by National Science Foundation with grant NSF-DMS-1504697.

\end{document}